\documentclass[11pt]{article}

\usepackage{latexsym}
\usepackage{amssymb}

\newcommand{\Z}{{\mathbb Z}} %
\newcommand{\R}{{\mathbb R}} %
\newcommand{\CC}{{\mathbb C}} %
\newcommand{\N}{{\mathbb N}} %
\newcommand{\PP}{{\mathbb P}} %
\newtheorem{thm}{Theorem}

\newcommand{\ld}{\lambda}
\newcommand{\ve}{\varepsilon}

\newcommand{\Lin}{\mathop{\rm Lin}\nolimits}
\newcommand{\img}{\mathop{\rm Im}\nolimits}

\newcommand{\mspan}{\mathop{\rm span}\nolimits}

\begin{document}

\begin{center}
{\Large\bf
A description of all solutions of the matrix Hamburger moment problem in a general case.}
\end{center}
\begin{center}
{\bf S.M. Zagorodnyuk}
\end{center}

\section{Introduction.}
The main aim of this investigation is to obtain a description of all solutions of the
matrix Hamburger moment problem. Recall that the matrix Hamburger moment problem consists of
finding a left-continuous non-decreasing matrix function $M(x) = ( m_{k,l}(x) )_{k,l=0}^{N-1}$
on $\R$, $M(-\infty)=0$, such that
\begin{equation}
\label{f1_1}
\int_\R x^n dM(x) = S_n,\qquad n\in\Z_+,
\end{equation}
where $\{ S_n \}_{n=0}^\infty$ is a given sequence of Hermitian $(N\times N)$ complex matrices, $N\in\N$.

Sequences $\{ S_n \}_{n=0}^\infty$ for which this problem has a solution are called {\it moment sequences}.
This problem was introduced in~1949 by M.G.~Krein~\cite{cit_100_K}. He described all solutions
in the case when the corresponding J-matrix defines a symmetric operator with maximal defect numbers.
This result appeared without proof in~\cite{cit_200_K}.
Using V.P.~Potapov's J-theory, in~1983 I.V.~Kovalishina described solutions of the matrix Hamburger
moment problem in the completely indeterminate case~\cite{cit_300_K} (The completely indeterminate case
meant that the limit radii of the matrix Weyl discs had full ranks).
Using properties of matrix orthogonal polynomials, in~2001 P.~Lopez-Rodriguez obtained a parameterization
of solutions in the completely indeterminate case~\cite{cit_400_L} (The completely indeterminate case
meant that the corresponding J-matrix generated a symmetric operator with maximal defect numbers).
In~2004, Yu.M.~Dyukarev introduced a notion of an abstract limit interpolation problem and
described solutions of the completely indeterminate limit interpolation problem~\cite{cit_500_D}.
As one of applications, he obtained a description of solutions of the matrix Hamburger
moment problem in the completely indeterminate case (This case
meant that the limit radii of the matrix Weyl discs had full ranks).

In the scalar case, a description of all solutions of the moment problem~(\ref{f1_1}) can be found, e.g.,
in~\cite{cit_600_Akh},\cite{cit_700_Ber} for the nondegenerate case, and in~\cite{cit_800_AK} for
the degenerate case.

Recall that the condition of solvability for the matrix Hamburger moment problem is that
for arbitrary complex vectors $\vec\xi_j = (\xi_{j,0},\xi_{j,1},\ldots,\xi_{j,N-1})$, $j=0,1,2,...$, it
holds (\cite[p. 52]{cit_100_K}):
\begin{equation}
\label{f1_2}
\sum_{j,k=0}^n \vec\xi_k^* S_{j+k} \vec\xi_j \geq 0,\qquad n=0,1,2,\ldots.
\end{equation}
Let us introduce the following matrices
\begin{equation}
\label{f1_3}
\Gamma_n = \left(
\begin{array}{cccc} S_0 & S_1 & \ldots & S_n\\
S_1 & S_2 & \ldots & S_{n+1}\\
\vdots & \vdots & \ddots & \vdots\\
S_n & S_{n+1} & \ldots & S_{2n}\end{array}
\right),\qquad n\in\Z_+.
\end{equation}
It is not hard to verify that condition~(\ref{f1_2}) is equivalent to the following inequalities
\begin{equation}
\label{f1_4}
\Gamma_n \geq 0,\qquad n\in\Z_+.
\end{equation}
In~1954, A.V.~Shtraus described all generalized resolvents of a densely defined symmetric operator
with an arbitrary deficiency index~\cite{cit_900_S}. In~1970, he described all generalized resolvents
for an arbitrary, not necessarily densely defined symmetric operator~\cite{cit_1000_S}.
We shall use these fundamental results to obtain a description of all solutions of the matrix
Hamburger moment problem in the case when condition~(\ref{f1_4}) is true.

We shall also study the truncated matrix Hamburger moment problem. The problem is
to find a left-continuous non-decreasing matrix function $M(x) = ( m_{k,l}(x) )_{k,l=0}^{N-1}$
on $\R$, $M(-\infty)=0$, such that
\begin{equation}
\label{f1_5}
\int_\R x^n dM(x) = S_n,\qquad n=0,1,\ldots,2d,
\end{equation}
where $\{ S_n \}_{n=0}^{2d}$ is a given sequence of Hermitian $(N\times N)$ complex matrices, $d\in\Z_+$, $N\in\N$.

The conditions of solvability of the moment problem~(\ref{f1_5}) were given by T.~Ando in~1970~\cite{cit_1100_A}.
The nondegenerate case of the truncated moment problem~(\ref{f1_5}) is the case when the following
condition takes place:
\begin{equation}
\label{f1_6}
\Gamma_{d}>0,
\end{equation}
where $\Gamma_{d}$ is defined as in~(\ref{f1_3}).
In~1968, V.G.~Ershov obtained a description of all solutions of the truncated matrix Hamburger moment
problem~(\ref{f1_5}) in the nondegenerate case, using an operator approach~\cite{cit_1200_E}.
In~1989, H.~Dym described all solutions of the moment problem~(\ref{f1_5}) in the nondegenerate case,
using the reproducing kernel Hilbert spaces approach~\cite{cit_1300_D}.
In~1997, V.M.~Adamyan and I.M.~Tkachenko obtained solutions of the truncated moment problem~(\ref{f1_5})
both in degenerate and nondegenerate cases, using an operator approach~\cite{cit_1400_AT}.
In~1998, G.-N. Chen and Y.-J. Hu obtained solutions of the truncated moment problem~(\ref{f1_5})
both in degenerate and nondegenerate cases, using a generalization of the Schur algorithm and
matrix  continued fractions~\cite{cit_1500_CH}.
We shall study the moment problem~(\ref{f1_5}) under the following conditions
\begin{equation}
\label{f1_7}
\Gamma_{d-1}>0,\quad \Gamma_{d}\geq 0,
\end{equation}
where $\Gamma_{d-1},\Gamma_d$ are defined as in~(\ref{f1_3}), $d\in\N$.
Using A.V.~Shtraus's results
we describe all solutions of the truncated moment problem~(\ref{f1_5}) under condition~(\ref{f1_7}).

Finally, we consider the scalar truncated Hamburger moment problem with even number of given moments.
The problem is to find a left-continuous non-decreasing function $\sigma(x)$
on $\R$, $\sigma(-\infty)=0$, such that
\begin{equation}
\label{f1_8}
\int_\R x^n d\sigma(x) = s_n,\qquad n=0,1,\ldots,2d+1,
\end{equation}
where $\{ s_n \}_{n=0}^{2d+1}$ is a given sequence of real numbers, $d\in\Z_+$.
Algebraic conditions of solvability of this moment problem were given in~\cite[Theorem 3.1]{cit_1600_CF}.
We shall give a simple condition of solvability for the truncated scalar Hamburger moment problem~(\ref{f1_7}).

For additional references on matrix Hamburger moment problems (including truncated) we refer to a historical
review in~\cite{cit_1700_D}.

\noindent
{\bf Notations.}  As usual, we denote by $\R, \CC, \N, \Z, \Z_+$
the sets of real, complex, positive integer, integer, non-negative integer numbers,
respectively. The space of $n$-dimensional complex vectors $a = (a_0,a_1,\ldots,a_{n-1})$, will be
denoted by $\CC^n$, $n\in\N$; $\CC_+ = \{ z\in\CC:\ \img z > 0\}$.
If $a\in \CC^n$ then $a^*$ means the complex conjugate vector.
By $\PP$ we denote a set of all complex polynomials and by $\PP_d$ we mean all complex polynomials with
degrees less or equal to $d$, $d\in\Z_+$, (including the zero polynomial).
Let $M(x)$ be a left-continuous non-decreasing matrix function $M(x) = ( m_{k,l}(x) )_{k,l=0}^{N-1}$
on $\R$, $M(-\infty)=0$, and $\tau_M (x) := \sum_{k=0}^{N-1} m_{k,k} (x)$;
$\Psi(x) = ( dm_{k,l}/ d\tau_M )_{k,l=0}^{N-1}$.  We denote by $L^2(M)$ a set (of classes of equivalence)
of vector functions $f: \R\rightarrow \CC^N$, $f = (f_0,f_1,\ldots,f_{N-1})$, such that (see, e.g.,~\cite{cit_1800_MM})
$$ \| f \|^2_{L^2(M)} := \int_\R  f(x) \Psi(x) f^*(x) d\tau_M (x) < \infty. $$
The space $L^2(M)$ is a Hilbert space with the scalar product
$$ ( f,g )_{L^2(M)} := \int_\R  f(x) \Psi(x) g^*(x) d\tau_M (x),\qquad f,g\in L^2(M). $$
By $l^2$ we denote a space of infinite complex vectors $u = (u_0,u_1,...)$, such that
$\| u \|_{l^2}^2 := \sum_{k=0}^\infty |u_k|^2 <\infty$.
The space $l^2$ is a Hilbert space with the scalar product $(u,v)_{l^2} = \sum_{k=0}^\infty u_k \overline{v_k}$,
$u,v\in l^2$. A set of elements $u = (u_0,u_1,...)$ from $l^2$, such that all but finite number
$u_k$ are zero will be denoted by $l^2_0$. Elements of $l^2_0$ are called finite vectors.

For a separable Hilbert space $H$ we denote by $(\cdot,\cdot)_H$ and $\| \cdot \|_H$ the scalar
product and the norm in $H$, respectively. The indices may be omitted in obvious cases.

\noindent
For a linear operator $A$ in $H$ we denote by $D(A)$ its domain, by $R(A)$ its range, and by
$A^*$ we denote its adjoint if it exists. If $A$ is bounded, then $\| A \|$ stands for its operator norm.
For a set of elements $\{ x_n \}_{n\in A}$ in $H$, we
denote by $\Lin\{ x_n \}_{n\in A}$ and $\mspan\{ x_n \}_{n\in A}$ the linear span and the closed
linear span (in the norm of $H$), respectively, where $A$ is an arbitrary set of indices.
For a set $M\subseteq H$ we denote by $\overline{M}$ the closure of $M$ with respect to the norm of $H$.
By $E_H$ we denote the identity operator in $H$, i.e. $E_H x = x$, $x\in H$.
If $H_1$ is a subspace of $H$, by $P_{H_1} = P_{H_1}^{H}$ we denote the operator of the orthogonal projection on $H_1$
in $H$.

\section{The matrix Hamburger moment problem: solvability and a description  of solutions.}
Recall that an infinite complex matrix $K= (K_{n,m})_{n,m=0}^\infty$ is called
a {\it positive definite kernel} if
\begin{equation}
\label{f2_1}
\sum_{n,m=0}^\infty K_{n,m} \xi_n \overline{\xi_m} \geq 0,
\end{equation}
for all finite vectors $(\xi_n)_{n=0}^\infty$ of complex numbers, see~\cite{cit_700_Ber}.
In other words, $K$ is a positive definite kernel if
\begin{equation}
\label{f2_2}
uKu^* = (uK, u )_{l^2} \geq 0,\quad u\in l^2_0,
\end{equation}
where $uK$ is defined by the usual matrix multiplication.

We shall use the following important fact (e.g., \cite[p.215]{cit_1900_K}):
\begin{thm}
\label{t2_1}
\noindent
a) Let $K = (K_{n,m})_{n,m=0}^\infty$ be a positive definite kernel.
Then there exist a separable Hilbert space $H$ with a scalar product $(\cdot,\cdot)$ and
a sequence $\{ x_n \}_{n=0}^\infty$ in $H$, such that
\begin{equation}
\label{f2_3}
K_{n,m} = (x_n,x_m),\qquad n,m\in \Z_+,
\end{equation}
and $\mspan\{ x_n \}_{n\in\Z_+} = H$.

\noindent
b) Let $R = (R_{n,m})_{n,m=0}^{r} \geq 0$ be a positive semi-definite complex $((r+1)\times(r+1))$ matrix, $r\in\Z_+$.
Then there exist a finite-dimensional Hilbert space $H_0$ with a scalar product $(\cdot,\cdot)_0$ and
a sequence $\{ y_n \}_{n=0}^r$ in $H_0$, such that
\begin{equation}
\label{f2_4}
R_{n,m} = (y_n,y_m),\qquad n,m=0,1,...,r,
\end{equation}
and $\mspan\{ y_n \}_{n=0}^r = H_0$.
\end{thm}
{\bf Proof. }
a) Consider an arbitrary infinite-dimensional linear vector space $V$ (for example a space of complex
sequences $(u_n)_{n\in\Z_+}$, $u_n\in\CC$). Let $X = \{ x_n \}_{n=0}^\infty$ be an arbitrary
infinite sequence of linear independent elements
in $V$. Let $L = \Lin\{ x_n \}_{n\in\Z_+}$ be the linear span of elements of $X$. Introduce the following functional:
\begin{equation}
\label{f2_5}
[x,y] = \sum_{n,m=0}^\infty K_{n,m} a_n\overline{b_m},
\end{equation}
for $x,y\in L$,
$$ x=\sum_{n=0}^\infty a_n x_n,\quad y=\sum_{m=0}^\infty b_m x_m,\quad a_n,b_m\in\CC. $$
The space $V$ with $[\cdot,\cdot]$ will be a quasi-Hilbert space. Factorizing and making the completion
we obtain the  required space $H$ (see~\cite[p. 10-11]{cit_700_Ber}).

\noindent
b) In this case we proceed in an analogous manner.
$\Box$

Consider the matrix Hamburger moment problem~(\ref{f1_1}).
If we choose an arbitrary element $f =(f_0,f_1,\ldots,f_{N-1})$, $f_k\in\PP$, $k=0,1,...,N-1$,
and calculate
$\int_\R f dM f^*$, one can easily deduce the necessity of conditions~(\ref{f1_2}),(\ref{f1_4})
for the solvability of the moment problem.

On the other hand, suppose that the moment problem~(\ref{f1_1}) is given and condition~(\ref{f1_4}) holds true.
Set
\begin{equation}
\label{f2_6}
\Gamma = (S_{k+l})_{k,l=0}^\infty = \left(
\begin{array}{ccccc} S_0 & S_1 & \ldots & S_n & \ldots\\
S_1 & S_2 & \ldots & S_{n+1} & \ldots\\
\vdots & \vdots & \ddots & \vdots & \ldots\\
S_n & S_{n+1} & \ldots & S_{2n} & \ldots\\
\vdots & \vdots & \vdots & \vdots & \ddots\end{array}
\right).
\end{equation}
Comparing relations~(\ref{f1_4}) and~(\ref{f2_2}) we conclude that the kernel
$\Gamma = (\Gamma_{n,m})_{n,m=0}^\infty$ is positive definite. Let
$$ S_n = (s_n^{k,l})_{k,l=0}^{N-1},\qquad n\in\Z_+. $$
Notice that
\begin{equation}
\label{f2_7}
\Gamma_{rN+j,tN+n} = s_{r+t}^{j,n},\qquad 0\leq j,n \leq N-1;\quad r,t\in\Z_+.
\end{equation}
From~(\ref{f2_7}) it follows that
\begin{equation}
\label{f2_8}
\Gamma_{a+N,b} = \Gamma_{a,b+N},\qquad a,b\in\Z_+.
\end{equation}
In fact, if $a=rN+j$, $b=tN+n$, $0\leq j,n \leq N-1$,  $r,t\in\Z_+$, we can write
$$ \Gamma_{a+N,b} = \Gamma_{(r+1)N+j,tN+n} = s_{r+t+1}^{j,n} = \Gamma_{rN+j,(t+1)N+n} = \Gamma_{a,b+N}. $$
By Theorem~\ref{t2_1} there exist a Hilbert space $H$ and a sequence $\{ x_n \}_{n=0}^\infty$ in $H$,
such that $\mspan\{ x_n \}_{n\in\Z_+} = H$, and
\begin{equation}
\label{f2_9}
(x_n,x_m)_H = \Gamma_{n,m},\qquad n,m\in\Z_+.
\end{equation}
Set $L := \Lin\{ x_n \}_{n\in\Z_+}$. Choose an arbitrary $x\in L$.
Let $x = \sum_{k=0}^\infty \alpha_k x_k$, $x = \sum_{k=0}^\infty \beta_k x_k$,
where $\alpha_k,\beta_k\in\CC$, and all but finite number of coefficients $\alpha_k$, $\beta_k$ are zero.
Using~(\ref{f2_9}),(\ref{f2_8}) we can write
$$ \left( \sum_{k=0}^\infty \alpha_k x_{k+N}, x_l \right) = \sum_{k=0}^\infty \alpha_k ( x_{k+N}, x_l ) =
\sum_{k=0}^\infty \alpha_k \Gamma_{k+N,l} = \sum_{k=0}^\infty \alpha_k \Gamma_{k,l+N} = $$
$$ = \sum_{k=0}^\infty \alpha_k ( x_{k}, x_{l+N} ) =
\left( \sum_{k=0}^\infty \alpha_k x_{k}, x_{l+N} \right) = (x,x_{l+N}),\qquad l\in\Z_+. $$
In an analogous manner we obtain that
$$ \left( \sum_{k=0}^\infty \beta_k x_{k+N}, x_l \right) = (x,x_{l+N}),\qquad l\in\Z_+, $$
and therefore
$$ \left( \sum_{k=0}^\infty \alpha_k x_{k+N}, x_l \right) =
\left( \sum_{k=0}^\infty \beta_k x_{k+N}, x_l \right),\qquad l\in\Z_+. $$
Since $\overline{L} = H$, we obtain that
\begin{equation}
\label{f2_10}
\sum_{k=0}^\infty \alpha_k x_{k+N} = \sum_{k=0}^\infty \beta_k x_{k+N}.
\end{equation}
Set
\begin{equation}
\label{f2_11}
A x = \sum_{k=0}^\infty \alpha_k x_{k+N},\qquad x\in L,\ x = \sum_{k=0}^\infty \alpha_k x_{k}.
\end{equation}
In particular, we have
\begin{equation}
\label{f2_12}
A x_k =  x_{k+N},\qquad k\in\Z_+.
\end{equation}
The above considerations show that this definition is correct.
Choose arbitrary $x,y\in L$, $x = \sum_{k=0}^\infty \alpha_k x_{k}$, $y = \sum_{n=0}^\infty \gamma_n x_{n}$,
and write
$$ (Ax,y) = \left( \sum_{k=0}^\infty \alpha_k x_{k+N},\sum_{n=0}^\infty \gamma_n x_{n} \right) =
\sum_{k,n=0}^\infty \alpha_k \overline{\gamma_n} (x_{k+N},x_n) =
\sum_{k,n=0}^\infty \alpha_k \overline{\gamma_n} (x_{k},x_{n+N}) = $$
$$ = \left( \sum_{k=0}^\infty \alpha_k x_{k},\sum_{n=0}^\infty \gamma_n x_{n+N} \right) =
(x,Ay). $$
Thus, the operator $A$ is a linear symmetric operator in $H$ with the domain $D(A)=L$.
Let $\widetilde A\supseteq A$ be an arbitrary self-adjoint extension of $A$ in a Hilbert space
$\widetilde H\supseteq H$, and $\{ \widetilde E_\ld \}_{\ld\in\R}$ be its left-continuous orthogonal
resolution of unity.
Choose an arbitrary $a\in\Z_+$, $a=rN + j$, $r\in\Z_+$, $0\leq j\leq N-1$. Notice that
$$ x_a = x_{rN+j} = A x_{(r-1)N+j} = ... = A^r x_j. $$
Then choose an arbitrary $b\in\Z_+$, $b=tN + n$, $t\in\Z_+$, $0\leq n\leq N-1$.
Using~(\ref{f2_7}) we can write
$$ s_{r+t}^{j,n} = \Gamma_{rN+j,tN+n} = ( x_{rN+j},x_{tN+n} )_H = (A^r x_j, A^t x_n)_H =
( \widetilde A^r x_j, \widetilde A^t x_n)_{\widetilde H} = $$
$$ = \left( \int_\R \ld^r d\widetilde E_\ld x_j, \int_\R \ld^t d\widetilde E_\ld x_n \right)_{\widetilde H} =
\int_\R \ld^{r+t} d (\widetilde E_\ld x_j, x_n)_{\widetilde H} =
\int_\R \ld^{r+t} d \left( P^{\widetilde H}_H \widetilde E_\ld x_j, x_n \right)_{H}. $$
From the latter relation we get
\begin{equation}
\label{f2_13}
S_{r+t} = \int_\R \ld^{r+t} d \widetilde M(\ld),\qquad r,t\in\Z_+,
\end{equation}
where $\widetilde M(\ld) := \left( \left( P^{\widetilde H}_H \widetilde E_\ld x_j,
x_n \right)_{H} \right)_{j,n=0}^{N-1}$.
If we set $t=0$ in relation~(\ref{f2_7}), we obtain that the matrix function $\widetilde M(\ld)$ is
a solution of the matrix Hamburger moment  problem~(\ref{f1_1}) (From the properties of the
orthogonal resolution of unity it easily follows that $\widetilde M (\ld)$ is left-continuous non-decreasing and
$\widetilde M(-\infty) = 0$).

Thus, we obtained another proof of the solvability criterion for the matrix Hamburger moment  problem~(\ref{f1_1}).

Let $\widehat A$ be an arbitrary self-adjoint extension of $A$ in a Hilbert space $\widehat H$.
Let $R_z(\widehat A)$ be the resolvent of $\widehat A$ and $\{ \widehat E_\lambda\}_{\ld\in\R}$
be an orthogonal left-continuous resolution of unity of $\widehat A$. Recall that the operator-valued function
$\mathbf R_z = P_H^{\widehat H} R_z(\widehat A)$ is called a generalized resolvent of $A$, $z\in\CC\backslash\R$.
The function
$\mathbf E_\lambda = P_H^{\widehat H} \widehat E_\ld$, $\ld\in\R$, is a spectral
function of a symmetric operator $A$.
There exists a one-to-one correspondence between generalized resolvents and spectral functions
established by the following relation (\cite{cit_2000_AG}):
\begin{equation}
\label{f2_14}
(\mathbf R_z f,g)_H = \int_\R \frac{1}{\ld - z} d( \mathbf E_\ld f,g)_H,\qquad f,g\in H,\ z\in\CC\backslash\R.
\end{equation}
Formula~(\ref{f2_13}) shows that spectral functions of $A$ produce solutions of the matrix Hamburger
moment problem~(\ref{f1_1}). Can an arbitrary solution of~(\ref{f1_1}) be produced in such a way?
Choose an arbitrary solution $\widehat M(x) = ( \widehat m_{k,l}(x) )_{k,l=0}^{N-1}$ of
the matrix Hamburger moment problem~(\ref{f1_1}). Consider the space $L^2(\widehat M)$ and
let $Q$ be the operator of multiplication by an independent variable in $L^2(\widehat M)$.
The operator $Q$ is self-adjoint and its resolution of unity is (see~\cite{cit_1800_MM})
\begin{equation}
\label{f2_14_1}
E_b - E_a = E([a,b)): h(x) \rightarrow \chi_{[a,b)}(x) h(x),
\end{equation}
where $\chi_{[a,b)}(x)$ is the characteristic function of an interval $[a,b)$, $-\infty\leq a<b\leq +\infty$.
Set $\vec e_k = (e_{k,0},e_{k,1},\ldots,e_{k,N-1})$, $e_{k,j}=\delta_{k,j}$, $0\leq j\leq N-1$,
for $k=0,1,\ldots N-1$.
A set of (classes of equivalence of) functions $f\in L^2(\widehat M)$ such that
(the corresponding class includes) $f=(f_0,f_1,\ldots, f_{N-1})$, $f\in\PP$, we denote
by $\PP^2(\widehat M)$ and call a set of vector polynomials in $L^2(\widehat M)$.
Set $L^2_0(\widehat M) = \overline{ \PP^2(\widehat M) }$.

For an arbitrary $f\in \PP^2(\widehat M)$ there exists a unique representation of the following form:
\begin{equation}
\label{f2_15}
f(x) = \sum_{k=0}^{N-1} \sum_{j=0}^\infty \alpha_{k,j} x^j \vec e_k,\quad (\alpha_{k,0},\alpha_{k,1},\ldots)\in l^2_0.
\end{equation}
Let $g\in \PP^2(\widehat M)$ have a representation
\begin{equation}
\label{f2_15_1}
g(x) = \sum_{l=0}^{N-1} \sum_{r=0}^\infty \beta_{l,r} x^r \vec e_l,\quad (\beta_{l,0},\beta_{l,1},\ldots)\in
l^2_0.
\end{equation}
We can write
$$ (f,g)_{L^2(\widehat M)} = \sum_{k,l=0}^{N-1} \sum_{j,r=0}^\infty \alpha_{k,j}\overline{\beta_{l,r}}
\int_\R x^{j+r} \vec e_k d\widehat M(x) \vec e_l^* = \sum_{k,l=0}^{N-1}
\sum_{j,r=0}^\infty \alpha_{k,j}\overline{\beta_{l,r}}
\int_\R x^{j+r} d\widehat m_{k,l}(x) = $$
\begin{equation}
\label{f2_16}
= \sum_{k,l=0}^{N-1} \sum_{j,r=0}^\infty \alpha_{k,j}\overline{\beta_{l,r}}
s_{j+r}^{k,l}.
\end{equation}
On the other hand, we can write
$$ \left( \sum_{j=0}^\infty \sum_{k=0}^{N-1} \alpha_{k,j} x_{jN+k},
\sum_{r=0}^\infty \sum_{l=0}^{N-1} \beta_{l,r} x_{rN+l} \right)_H =
\sum_{k,l=0}^{N-1} \sum_{j,r=0}^\infty \alpha_{k,j}\overline{\beta_{l,r}}
(x_{jN+k}, x_{rN+l})_H  = $$
\begin{equation}
\label{f2_17}
= \sum_{k,l=0}^{N-1} \sum_{j,r=0}^\infty \alpha_{k,j}\overline{\beta_{l,r}}
\Gamma_{jN+k,rN+l}
= \sum_{k,l=0}^{N-1} \sum_{j,r=0}^\infty \alpha_{k,j}\overline{\beta_{l,r}}
s_{j+r}^{k,l}.
\end{equation}
From relations~(\ref{f2_16}),(\ref{f2_17}) it follows that
\begin{equation}
\label{f2_18}
(f,g)_{L^2(\widehat M)} = \left( \sum_{j=0}^\infty \sum_{k=0}^{N-1} \alpha_{k,j} x_{jN+k},
\sum_{r=0}^\infty \sum_{l=0}^{N-1} \beta_{l,r} x_{rN+l} \right)_H.
\end{equation}
Set
\begin{equation}
\label{f2_19}
Vf = \sum_{j=0}^\infty \sum_{k=0}^{N-1} \alpha_{k,j} x_{jN+k},
\end{equation}
for $f(x)\in \PP^2(\widehat M)$, $f(x) = \sum_{k=0}^{N-1} \sum_{j=0}^\infty \alpha_{k,j} x^j \vec e_k$,
$(\alpha_{k,0},\alpha_{k,1},\ldots)\in l^2_0$.

\noindent
If $f$, $g$ have representations~(\ref{f2_15}),(\ref{f2_15_1}), and $\| f-g \|_{L^2(\widehat M)} = 0$, then
from~(\ref{f2_18}) it follows that
$$ \| Vf - Vg \|_H^2 = (V(f-g),V(f-g))_H = ( f-g,f-g )_{L^2(\widehat M)} = \| f-g\|_{L^2(\widehat M)}^2 = 0.$$
Thus, $V$ is a correctly defined operator from $\PP^2(\widehat M)$ to $H$.
Relation~(\ref{f2_18}) shows that $V$ is an isometric transformation from $\PP^2(\widehat M)$ onto $L$.
By continuity we extend it to an isometric transformation from $L^2_0(\widehat M)$ onto $H$.
In particular, we note that
\begin{equation}
\label{f2_20}
V x^j \vec e_k = x_{jN+k},\qquad j\in\Z_+;\quad 0\leq k\leq N-1.
\end{equation}
Set $L^2_1 (\widehat M) := L^2(\widehat M)\ominus L^2_0 (\widehat M)$, and
$U := V\oplus E_{L^2_1 (\widehat M)}$. The operator $U$ is
an isometric transformation from $L^2(\widehat M)$ onto $H\oplus L^2_1 (\widehat M)=:\widehat H$.
Set
$$ \widehat A := UQU^{-1}. $$
The operator $\widehat A$ is a self-adjoint operator in $\widehat H$. Let $\{ \widehat E_\ld \}_{\ld\in\R}$
be its left-continuous orthogonal resolution of unity.
Notice that
$$ UQU^{-1} x_{jN+k} = VQV^{-1} x_{jN+k} = VQ x^j \vec e_k = V x^{j+1} \vec e_k =
x_{(j+1)N+k} = x_{jN+k+N} = $$
$$ = Ax_{jN+k},\qquad j\in\Z_+;\quad 0\leq k\leq N-1. $$
By linearity we get
$$ UQU^{-1} x = Ax,\qquad x\in L = D(A), $$
and therefore $\widehat A\supseteq A$.
Choose an arbitrary $z\in\CC\backslash\R$ and write
$$ \int_\R \frac{1}{\ld - z} d( \widehat E_\ld x_k, x_j)_{\widehat H} =
\left( \int_\R \frac{1}{\ld - z} d\widehat E_\ld x_k, x_j \right)_{\widehat H} =
\left( U^{-1} \int_\R \frac{1}{\ld - z} d\widehat E_\ld x_k, U^{-1} x_j \right)_{L^2(\widehat M)} = $$
$$ = \left( \int_\R \frac{1}{\ld - z} d U^{-1} \widehat E_\ld U \vec e_k, \vec e_j \right)_{L^2(\widehat M)} =
\left( \int_\R \frac{1}{\ld - z} d E_{\ld} \vec e_k, \vec e_j \right)_{L^2(\widehat M)} = $$
\begin{equation}
\label{f2_21}
= \int_\R \frac{1}{\ld - z} d(E_{\ld} \vec e_k, \vec e_j)_{L^2(\widehat M)},\qquad 0\leq k,j\leq N-1.
\end{equation}
Using~(\ref{f2_14_1}) we can write
$$ (E_{\ld} \vec e_k, \vec e_j)_{L^2(\widehat M)} = \widehat m_{k,j}(\ld), $$
and therefore
\begin{equation}
\label{f2_22}
\int_\R \frac{1}{\ld - z} d( P^{\widehat H}_H \widehat E_\ld x_k, x_j)_H =
\int_\R \frac{1}{\ld - z} d\widehat m_{k,j}(\ld),\qquad 0\leq k,j\leq N-1.
\end{equation}
By the Stieltjes-Perron inversion  formula (see, e.g., \cite{cit_600_Akh}) we conclude that
\begin{equation}
\label{f2_23}
\widehat m_{k,j} (\ld) = ( P^{\widehat H}_H \widehat E_\ld x_k, x_j)_H.
\end{equation}
Consequently, an answer on the above question is affirmative.

Let us show that the deficiency index of $A$ is equal to $(m,n)$, $0\leq m,n\leq N$.
Choose an arbitrary $u\in L$, $u = \sum_{k=0}^\infty c_k x_k$, $c_k\in\CC$. Suppose that
$c_k = 0$, $k\geq N+R+1$, for some $R\in\Z_+$. Consider the following system of linear equations:
\begin{equation}
\label{f2_24}
-z d_k = c_k,\qquad  k=0,1,...,N-1;
\end{equation}
\begin{equation}
\label{f2_25}
d_{k-N} - z d_k = c_k,\qquad  k=N,N+1,N+2,...;
\end{equation}
where $\{ d_k \}_{k\in\Z_+}$ are unknown complex numbers, $z\in\CC\backslash\R$ is a fixed parameter.
Set
$$ d_k = 0,\qquad k\geq R+1; $$
\begin{equation}
\label{f2_26}
d_{j} = c_{N+j} + z d_{N+j},\qquad j=R,R-1,R-2,...,0.
\end{equation}
For such numbers $\{ d_k \}_{k\in\Z_+}$, all equations in~(\ref{f2_25}) are satisfied.
Only equations~(\ref{f2_24}) are not satisfied. Set
$v = \sum_{k=0}^\infty d_k x_k$, $v\in L$.
Notice that
$$ (A-zE_H) v = \sum_{k=0}^\infty (d_{k-N} - z d_k) x_k, $$
where $d_{-1}=d_{-2}=...=d_{-N}=0$.
By the construction of $d_k$ we have
$$ (A-zE_H) v - u = \sum_{k=0}^\infty (d_{k-N} - z d_k - c_k) x_k =
\sum_{k=0}^{N-1} (-zd_k - c_k) x_k; $$
\begin{equation}
\label{f2_27}
u = (A-zE_H) v + \sum_{k=0}^{N-1} (zd_k + c_k) x_k,\qquad u\in L.
\end{equation}
Set $H_z := \overline{(A-zE_H) L} = (\overline{A} - zE_H) D(\overline{A})$, and
\begin{equation}
\label{f2_28}
y_k := x_k - P^H_{H_z} x_k,\qquad k=0,1,...,N-1.
\end{equation}
Set $H_0 := \mspan\{ y_k \}_{k=0}^{N-1}$. Notice that the dimension of $H_0$ is less or equal to
$N$, and $H_0\perp H_z$.
From~(\ref{f2_27}) it follows that $u\in L$ can be represented in the following form:
\begin{equation}
\label{f2_29}
u = u_1 + u_2,\qquad u_1\in H_z,\quad u_2\in H_0.
\end{equation}
Therefore we get $L\subseteq H_z\oplus H_0$; $H\subseteq H_z\oplus H_0$, and finally
$H=H_z\oplus H_0$. Thus, $H_0$ is the corresponding defect subspace.
So, the defect numbers of $A$ are less or equal to $N$.

\begin{thm}
\label{t2_2}
Let a matrix Hamburger moment problem~(\ref{f1_1}) be given and
condition~(\ref{f1_4}) is true. Let an operator $A$ be constructed for the
moment problem as in~(\ref{f2_11}).
All solutions of the moment problem have the following form
\begin{equation}
\label{f2_30}
M(\ld) = (m_{k,j} (\ld))_{k,j=0}^{N-1},\quad
m_{k,j} (\ld) = ( \mathbf E_\ld x_k, x_j)_H,
\end{equation}
where $\mathbf E_\ld$ is a spectral function of the operator $A$.
Moreover, the correspondence between all spectral functions of $A$ and all solutions
of the moment problem is one-to-one.
\end{thm}
{\bf Proof. }
It remains to prove that different spectral functions of the operator $A$ produce different
solutions of the moment problem~(\ref{f1_1}).
Suppose to the contrary that two different spectral functions produce the same solution of
the moment problem. That means that
there exist two self-adjoint extensions
$A_j\supseteq A$, in Hilbert spaces $H_j\supseteq H$, such that
\begin{equation}
\label{f2_30_1}
P_{H}^{H_1} E_{1,\ld} \not= P_{H}^{H_2} E_{2,\ld},
\end{equation}
\begin{equation}
\label{f2_31}
(P_{H}^{H_1} E_{1,\ld} x_k,x_j)_H = (P_{H}^{H_2} E_{2,\ld} x_k,x_j)_H,\qquad 0\leq k,j\leq N-1,\quad \ld\in\R,
\end{equation}
where $\{ E_{n,\ld} \}_{\ld\in\R}$ are orthogonal left-continuous resolutions of unity of
operators $A_n$, $n=1,2$.
Set $L_N := \Lin\{ x_k \}_{k=0,N-1}$. By linearity we get
\begin{equation}
\label{f2_32}
(P_{H}^{H_1} E_{1,\ld} x,y)_H = (P_{H}^{H_2} E_{2,\ld} x,y)_H,\qquad x,y\in L_N,\quad \ld\in\R.
\end{equation}
Denote by $R_{n,\ld}$ the resolvent of $A_n$, and set $\mathbf R_{n,\ld} := P_{H}^{H_n} R_{n,\ld}$, $n=1,2$.
From~(\ref{f2_32}),(\ref{f2_14}) it follows that
\begin{equation}
\label{f2_33}
(\mathbf R_{1,\ld} x,y)_H = (\mathbf R_{2,\ld} x,y)_H,\qquad x,y\in L_N,\quad \ld\in\CC\backslash\R.
\end{equation}
Choose an arbitrary $z\in\CC\backslash\R$ and consider the space $H_z$ defined as above.
Since
$$ R_{j,z} (A-zE_H) x = (A_j - z E_{H_j} )^{-1} (A_j - z E_{H_j}) x = x,\qquad x\in L=D(A),$$
we get
\begin{equation}
\label{f2_34}
R_{1,z} u = R_{2,z} u \in H,\qquad u\in H_z;
\end{equation}
\begin{equation}
\label{f2_35}
\mathbf R_{1,z} u = \mathbf R_{2,z} u,\qquad u\in H_z,\ z\in\CC\backslash\R.
\end{equation}
We can write
\begin{equation}
\label{f2_36}
(\mathbf R_{n,z} x, u)_H = (R_{n,z} x, u)_{H_n} = ( x, R_{n,\overline{z}}u)_{H_n} =
( x, \mathbf R_{n,\overline{z}} u)_H,\qquad x\in L_N,\ u\in H_{\overline z},\ n=1,2,
\end{equation}
and therefore we get
\begin{equation}
\label{f2_37}
(\mathbf R_{1,z} x,u)_H = (\mathbf R_{2,z} x,u)_H,\qquad x\in L_N,\ u\in H_{\overline z}.
\end{equation}
By~(\ref{f2_27}) an arbitrary element $y\in L$ can be represented as $y=y_{ \overline{z} } + y'$,
$y_{ \overline{z} }\in H_{ \overline{z} }$, $y'\in L_N$.
Using~(\ref{f2_33}) and~(\ref{f2_37})  we get
$$ (\mathbf R_{1,z} x,y)_H = (\mathbf R_{1,z} x, y_{ \overline{z} } + y')_H =
(\mathbf R_{2,z} x, y_{ \overline{z} } + y')_H = (\mathbf R_{2,z} x,y)_H,\qquad x\in L_N,\ y\in L. $$
Since $\overline{L}=H$, we obtain
\begin{equation}
\label{f2_38}
\mathbf R_{1,z} x = \mathbf R_{2,z} x,\qquad x\in L_N,\ z\in\CC\backslash\R.
\end{equation}
For an arbitrary $x\in L$, $x=x_z + x'$, $x_z\in H_z$, $x'\in L_N$, using
relations~(\ref{f2_35}),(\ref{f2_38}) we obtain
\begin{equation}
\label{f2_39}
\mathbf R_{1,z} x = \mathbf R_{1,z} (x_z + x') =
\mathbf R_{2,z} (x_z + x') = \mathbf R_{2,z} x,\qquad x\in L,\ z\in\CC\backslash\R,
\end{equation}
and
\begin{equation}
\label{f2_40}
\mathbf R_{1,z} x = \mathbf R_{2,z} x,\qquad x\in H,\ z\in\CC\backslash\R.
\end{equation}
By~(\ref{f2_14}) that means that the spectral functions coincide and we obtain a
contradiction.
$\Box$

Recall some known facts from~\cite{cit_900_S} which we shall need here.
Let $B$ be a closed symmetric operator in a Hilbert space $H$, with the domain $D(B)$,
$\overline{D(B)} = H$.  Set $\Delta_B(\ld) = (B-\ld E_H) D(B)$,
and $N_\ld = N_\ld(B) = H\ominus \Delta_B(\ld)$, $\ld\in\CC\backslash\R$.

Consider an
arbitrary bounded linear operator $C$, which maps $N_i$ into $N_{-i}$.
For
\begin{equation}
\label{f2_41}
g = f + C\psi - \psi,\qquad f\in D(B),\ \psi\in N_i,
\end{equation}
we set
\begin{equation}
\label{f2_42}
B_C g = Bf + i C \psi + i \psi.
\end{equation}
Since an intersection of $D(A)$, $N_i$ and $N_{-i}$ consists only of the zero element,
this definition is correct.
Notice that $B_C$ is a part of the operator $B^*$.
The operator $B_C$ is called a {\it quasiself-adjoint extension of the operator $B$, defined by
the operator $C$}.

The following theorem is true, see~\cite[Theorem 7]{cit_900_S}:
\begin{thm}
\label{t2_3}
Let $B$ be a closed symmetric operator in a Hilbert space $H$ with the domain $D(B)$,
$\overline{D(B)} = H$.
All generalized resolvents of the operator $B$ have the following form:
\begin{equation}
\label{f2_43}
\mathbf R_\ld = \left\{ \begin{array}{cc} (B_{F(\ld)} - \ld E_H)^{-1}, & \img\ld > 0\\
(B_{F^*(\overline{\ld}) } - \ld E_H)^{-1}, & \img\ld < 0 \end{array}\right.,
\end{equation}
where $F(\ld)$ is an analytic in $\CC_+$ operator-valued function, which values are contractions
which map $N_i(B)$ into $N_{-i}(B)$ ($\| F(\ld) \|\leq 1$),
and $B_{F(\ld)}$ is the quasiself-adjoint extension of $B$ defined by $F(\ld)$.

On the other hand, for any operator function $F(\ld)$ having the above properties there corresponds by
relation~(\ref{f2_43}) a generalized resolvent of $B$.
\end{thm}

By virtue of Theorems~\ref{t2_2} and~\ref{t2_3} we get a description of all solutions of the
matrix Hamburger moment problem~(\ref{f1_1}).
\begin{thm}
\label{t2_4}
Let a matrix Hamburger moment problem~(\ref{f1_1}) be given and
condition~(\ref{f1_4}) is true. Let an operator $A$ be constructed for the
moment problem as in~(\ref{f2_11}).
All solutions of the moment problem have the following form
\begin{equation}
\label{f2_44}
M(x) = (m_{k,j} (x))_{k,j=0}^{N-1},
\end{equation}
where $m_{k,j}$ satisfy the following relation
\begin{equation}
\label{f2_45}
\int_R \frac{1}{x-\ld} d m_{k,j} (x) = ( (A_{F(\ld)} - \ld E_H)^{-1} x_k, x_j)_H,\qquad \ld\in\CC_+,
\end{equation}
where $F(\ld)$ is an analytic in $\CC_+$ operator-valued function, which values are contractions
which map $N_i(\overline{A})$ into $N_{-i}(\overline{A})$ ($\| F(\ld) \|\leq 1$), and $A_{F(\ld)}$ is the
quasiself-adjoint extension of $\overline{A}$ defined by $F(\ld)$.

On the other hand, to any operator function $F(\ld)$ having the above properties there corresponds by
relation~(\ref{f2_45}) a solution of the matrix Hamburger moment problem.
Moreover, the correspondence between all operator functions having the above properties and all solutions
of the moment problem, established by relation~(\ref{f2_45}), is one-to-one.
\end{thm}
{\bf Proof. } It remains to check the last statement of the theorem.
Note that different functions $F_1(\ld)$, $F_2(\ld)$, with the above properties
generate different generalized resolvents $\mathbf R_1 (\ld)$,
$\mathbf R_2 (\ld)$ of $A$ (see~\cite[Remark 2, p. 85]{cit_900_S}).
Let $\mathbf E_1 (\ld)$, $\mathbf E_2 (\ld)$, be the corresponding spectral functions of $A$.
Suppose to the contrary that functions $F_1(\ld)$, $F_2(\ld)$, correspond to the same solution
$M(x) = (m_{k,j} (x))_{k,j=0}^{N-1}$ of the moment problem. By~(\ref{f2_45}) this means that
\begin{equation}
\label{f2_46}
\int_R \frac{1}{x-\ld} d m_{k,j} (x) = ( \mathbf R_1 (\ld) x_k, x_j )_H = ( \mathbf R_2 (\ld) x_k, x_j )_H.
\end{equation}
By the Stiltjes-Perron inversion formula we get
\begin{equation}
\label{f2_47}
m_{k,j} (x) = ( \mathbf E_1 (\ld) x_k, x_j )_H = ( \mathbf E_2 (\ld) x_k, x_j )_H.
\end{equation}
We obtain that different spectral functions of $A$ generate the same solution of the moment problem.
This contradicts to Theorem~\ref{t2_2}.
$\Box$

\section{The truncated matrix Hamburger moment problem.}
Let a moment problem~(\ref{f1_5}) be given with $d\in\N$, and condition~(\ref{f1_7}) is true.
Let $\Gamma_d = (\gamma_{n,m}^d)_{n,m=0}^{dN +N-1}$.
By~Theorem~\ref{t2_1} there exist a finite-dimensional Hilbert space $H$ and
a sequence $\{ x_n \}_{n=0}^{dN +N-1}$ in $H$, such that
\begin{equation}
\label{f3_1}
\gamma_{n,m}^d = (x_n,x_m),\qquad n,m=0,1,...,dN +N-1,
\end{equation}
and $\mspan\{ x_n \}_{n=0}^{dN +N-1} = H$.
Notice that
\begin{equation}
\label{f3_2}
\gamma_{rN+j,tN+n}^d = s_{r+t}^{j,n},\qquad 0\leq j,n \leq N-1;\quad 0\leq r,t\leq d.
\end{equation}
From~(\ref{f3_2}) it follows that
\begin{equation}
\label{f3_3}
\gamma_{a+N,b}^d = \gamma_{a,b+N}^d,\qquad a=rN+j,\ b=tN+n,\ 0\leq j,n \leq N-1;\quad 0\leq r,t\leq d-1.
\end{equation}
In fact, we can write
$$ \gamma_{a+N,b}^d = \gamma_{(r+1)N+j,tN+n}^d = s_{r+t+1}^{j,n} = \gamma_{rN+j,(t+1)N+n}^d =
\gamma_{a,b+N}^d. $$
The first relation in~(\ref{f1_7}) means that the Gram matrix of elements $\{ x_n \}_{n=0}^{dN-1}$ is
positive. Therefore these elements are linear independent. Denote $L_a = \Lin\{ x_n \}_{n=0}^{dN-1}$,
$L = \Lin\{ x_n \}_{n=0}^{dN+N-1}$.
Set
\begin{equation}
\label{f3_4}
A x = \sum_{k=0}^{dN-1} \alpha_k x_{k+N},\qquad x\in L_a,\ x = \sum_{k=0}^{dN-1} \alpha_k x_{k}.
\end{equation}
In particular, we have
\begin{equation}
\label{f3_5}
A x_k =  x_{k+N},\qquad 0\leq k\leq dN-1.
\end{equation}
Choose arbitrary $x,y\in L_a$, $x = \sum_{k=0}^{dN-1} \alpha_k x_{k}$, $y = \sum_{n=0}^{dN-1} \gamma_n x_{n}$,
and write
$$ (Ax,y) = \left( \sum_{k=0}^{dN-1} \alpha_k x_{k+N},\sum_{n=0}^{dN-1} \gamma_n x_{n} \right) =
\sum_{k,n=0}^{dN-1} \alpha_k \overline{\gamma_n} (x_{k+N},x_n) =
\sum_{k,n=0}^{dN-1} \alpha_k \overline{\gamma_n} (x_{k},x_{n+N}) = $$
$$ = \left( \sum_{k=0}^{dN-1} \alpha_k x_{k},\sum_{n=0}^{dN-1} \gamma_n x_{n+N} \right) =
(x,Ay). $$
Thus, an operator $A$ is a linear symmetric operator in $H$ with the domain $D(A)=L_a$.
It is not necessary that $A$ is densely defined.

Let $\widetilde A\supseteq A$ be an arbitrary self-adjoint extension of $A$ in a Hilbert space
$\widetilde H\supseteq H$, and $\{ \widetilde E_\ld \}_{\ld\in\R}$ be its left-continuous orthogonal
resolution of unity. Existence of a self-adjoint extension of a non-densely defined symmetric operator
was established by M.A.~Krasnoselskiy (e.g.~\cite{cit_900_S}).
Choose an arbitrary $a$, $0\leq a\leq dN+N-1$, $a=rN + j$, $0\leq r\leq d$, $0\leq j\leq N-1$. Notice that
$$ x_a = x_{rN+j} = A x_{(r-1)N+j} = ... = A^r x_j. $$
Then choose an arbitrary $b$, $0\leq b\leq dN+N-1$, $b=tN + n$, $0\leq t\leq d$, $0\leq n\leq N-1$.
Using~(\ref{f3_1}) we can write
$$ s_{r+t}^{j,n} = \gamma_{rN+j,tN+n}^d = ( x_{rN+j},x_{tN+n} )_H = (A^r x_j, A^t x_n)_H =
( \widetilde A^r x_j, \widetilde A^t x_n)_{\widetilde H} = $$
$$ = \left( \int_\R \ld^r d\widetilde E_\ld x_j, \int_\R \ld^t d\widetilde E_\ld x_n \right)_{\widetilde H} =
\int_\R \ld^{r+t} d (\widetilde E_\ld x_j, x_n)_{\widetilde H} =
\int_\R \ld^{r+t} d \left( P^{\widetilde H}_H \widetilde E_\ld x_j, x_n \right)_{H}. $$
From the last relation we obtain
\begin{equation}
\label{f3_6}
S_{r+t} = \int_\R \ld^{r+t} d \widetilde M(\ld),\qquad 0\leq r,t\leq d,
\end{equation}
where $\widetilde M(\ld) := \left( \left( P^{\widetilde H}_H \widetilde E_\ld x_j,
x_n \right)_{H} \right)_{j,n=0}^{N-1}$.
From~relation~(\ref{f2_7}) we derive that the matrix function $\widetilde M(\ld)$ is
a solution of the matrix Hamburger moment  problem~(\ref{f1_5}) (Properties of the
orthogonal resolution of unity provide that $\widetilde M (\ld)$ is left-continuous non-decreasing and
$\widetilde M(-\infty) = 0$).

On the other hand, choose an arbitrary solution $\widehat M(x) = ( \widehat m_{k,l}(x) )_{k,l=0}^{N-1}$ of
the truncated matrix Hamburger moment problem~(\ref{f1_5}). Consider the space $L^2(\widehat M)$ and
let $Q$ be the operator of multiplication by an independent variable in $L^2(\widehat M)$.
The operator $Q$ is self-adjoint and its resolution of unity is given by~(\ref{f2_14_1}).

Let $\vec e_k$, $k=0,1,\ldots N-1$, be defined as after~(\ref{f2_14_1}).
A set of (classes of equivalence of) functions $f\in L^2(\widehat M)$ such that
(the corresponding class includes) $f=(f_0,f_1,\ldots, f_{N-1})$, $f\in\PP_d$, we denote
by $\PP^2_d(\widehat M)$ and call a set of vector polynomials in $L^2(\widehat M)$ of degree less or equal
to $d$.
Set $L^2_{d,0}(\widehat M) = \overline{ \PP^2_d(\widehat M) }$.

For an arbitrary $f\in \PP^2_d(\widehat M)$ there exists a unique representation of the following form:
\begin{equation}
\label{f3_7}
f(x) = \sum_{k=0}^{N-1} \sum_{j=0}^d \alpha_{k,j} x^j \vec e_k,\quad \alpha_{k,j}\in\CC.
\end{equation}
Let $g\in \PP^2_d(\widehat M)$ has a representation
$$ g(x) = \sum_{l=0}^{N-1} \sum_{r=0}^d \beta_{l,r} x^r \vec e_l,\quad \beta_{l,r}\in\CC. $$
As it was done in the  case of the full matrix Hamburger moment problem after~(\ref{f2_14_1}), we
obtain that
\begin{equation}
\label{f3_10}
(f,g)_{L^2(\widehat M)} = \left( \sum_{j=0}^d \sum_{k=0}^{N-1} \alpha_{k,j} x_{jN+k},
\sum_{r=0}^d \sum_{l=0}^{N-1} \beta_{l,r} x_{rN+l} \right)_H.
\end{equation}
Set
\begin{equation}
\label{f3_11}
Vf = \sum_{j=0}^d \sum_{k=0}^{N-1} \alpha_{k,j} x_{jN+k},
\end{equation}
for $f(x) = \sum_{k=0}^{N-1} \sum_{j=0}^d \alpha_{k,j} x^j \vec e_k$,
$(\alpha_{k,0},\alpha_{k,1},\ldots)\in l^2_0$.
From relation~(\ref{f3_10}) it easily follows that $V$ is a correctly defined operator from $\PP^2_d(\widehat M)$ to $H$.
Relation~(\ref{f3_10}) shows that $V$ is an isometric transformation from $\PP^2_d(\widehat M)$ onto $L$.
By continuity we extend it to an isometric transformation from $L^2_{d,0}(\widehat M)$ onto $H$.
In particular, we note that
\begin{equation}
\label{f3_12}
V x^j \vec e_k = x_{jN+k},\qquad 0\leq j\leq d;\quad 0\leq k\leq N-1.
\end{equation}
Set $L^2_{d,1} (\widehat M) := L^2(\widehat M)\ominus L^2_{d,0} (\widehat M)$, and
$U := V\oplus E_{L^2_{d,1} (\widehat M)}$. The operator $U$ is
an isometric transformation from $L^2(\widehat M)$ onto $H\oplus L^2_{d,1} (\widehat M)=:\widehat H$.
Set
$$ \widehat A := UQU^{-1}. $$
The operator $\widehat A$ is a self-adjoint operator in $\widehat H$. Let $\{ \widehat E_\ld \}_{\ld\in\R}$
be its left-continuous orthogonal resolution of unity.
Notice that
$$ UQU^{-1} x_{jN+k} = VQV^{-1} x_{jN+k} = VQ x^j \vec e_k = V x^{j+1} \vec e_k =
x_{(j+1)N+k} = x_{jN+k+N} = $$
$$ = Ax_{jN+k},\qquad 0\leq j\leq d-1;\quad 0\leq k\leq N-1. $$
By linearity we get
$$ UQU^{-1} x = Ax,\qquad x\in L_a = D(A), $$
and therefore $\widehat A\supseteq A$.
Choose an arbitrary $z\in\CC\backslash\R$ and write
$$ \int_\R \frac{1}{\ld - z} d( \widehat E_\ld x_k, x_j)_{\widehat H} =
\left( \int_\R \frac{1}{\ld - z} d\widehat E_\ld x_k, x_j \right)_{\widehat H} =
\left( U^{-1} \int_\R \frac{1}{\ld - z} d\widehat E_\ld x_k, U^{-1} x_j \right)_{L^2(\widehat M)} = $$
$$ = \left( \int_\R \frac{1}{\ld - z} d U^{-1} \widehat E_\ld U \vec e_k, \vec e_j \right)_{L^2(\widehat M)} =
\left( \int_\R \frac{1}{\ld - z} d E_{\ld} \vec e_k, \vec e_j \right)_{L^2(\widehat M)} = $$
\begin{equation}
\label{f3_13}
= \int_\R \frac{1}{\ld - z} d(E_{\ld} \vec e_k, \vec e_j)_{L^2(\widehat M)},\qquad 0\leq k,j\leq N-1.
\end{equation}
Using~(\ref{f2_14_1}) we can write
$$ (E_{\ld} \vec e_k, \vec e_j)_{L^2(\widehat M)} = \widehat m_{k,j}(\ld), $$
and therefore
\begin{equation}
\label{f3_14}
\int_\R \frac{1}{\ld - z} d( P^{\widehat H}_H \widehat E_\ld x_k, x_j)_H =
\int_\R \frac{1}{\ld - z} d\widehat m_{k,j}(\ld),\qquad 0\leq k,j\leq N-1.
\end{equation}
By the Stieltjes-Perron inversion  formula we conclude that
\begin{equation}
\label{f3_15}
\widehat m_{k,j} (\ld) = ( P^{\widehat H}_H \widehat E_\ld x_k, x_j)_H.
\end{equation}
Consequently, all solutions of the truncated moment problem are generated by
spectral functions of $A$. For the definitions of a spectral function and a generalized resolvent
for a non-densely defined symmetric operator we refer to~\cite{cit_900_S}.

Let us show that the deficiency index of $A$ is equal to $(m,n)$, $0\leq m,n\leq N$.
Choose an arbitrary $u\in L$, $u = \sum_{k=0}^{dN+N-1} c_k x_k$, $c_k\in\CC$.
Consider the following system of linear equations:
\begin{equation}
\label{f3_16}
-z d_k = c_k,\qquad  k=0,1,...,N-1;
\end{equation}
\begin{equation}
\label{f3_17}
d_{k-N} - z d_k = c_k,\qquad  k=N,N+1,\ldots,dN+N-1;
\end{equation}
where $\{ d_k \}_{k=0}^{dN+N-1}$ are unknown complex numbers, $z\in\CC\backslash\R$ is a fixed parameter.
Set
$$ d_k = 0,\qquad k=dN,dN+1,...,dN+N-1; $$
\begin{equation}
\label{f3_18}
d_{k-N} = z d_k + c_k,\qquad  k=dN+N-1,dN+N-2,...,N;
\end{equation}
For such numbers $\{ d_k \}_{k\in\Z_+}$, all equations in~(\ref{f3_17}) are satisfied.
Equations~(\ref{f3_16}) are not necessarily satisfied. Set
$v = \sum_{k=0}^{dN+N-1} d_k x_k = \sum_{k=0}^{dN-1} d_k x_k$. Notice that $v\in L_a = D(A)$.
We can write
$$ (A-zE_H) v = \sum_{k=0}^{dN+N-1} (d_{k-N} - z d_k) x_k, $$
where $d_{-1}=d_{-2}=...=d_{-N}=0$.
By the construction of $d_k$ we have
$$ (A-zE_H) v - u = \sum_{k=0}^{dN+N-1} (d_{k-N} - z d_k - c_k) x_k =
\sum_{k=0}^{N-1} (-zd_k - c_k) x_k; $$
\begin{equation}
\label{f3_19}
u = (A-zE_H) v + \sum_{k=0}^{N-1} (zd_k + c_k) x_k,\qquad u\in L.
\end{equation}
Set $H_z := \overline{(A-zE_H) L} = (\overline{A} - zE_H) D(\overline{A})$.
Repeating arguments after relation~(\ref{f2_27}) we obtain that the
defect numbers of $A$ are less or equal to $N$.

\begin{thm}
\label{t3_1}
Let a truncated matrix Hamburger moment problem~(\ref{f1_5}) with $d\in\N$ be given and
conditions~(\ref{f1_7}) are true. Let the operator $A$ be constructed for the
moment problem as in~(\ref{f3_4}).
All solutions of the moment problem have the following form
\begin{equation}
\label{f3_22}
M(\ld) = (m_{k,j} (\ld))_{k,j=0}^{N-1},\quad
m_{k,j} (\ld) = ( \mathbf E_\ld x_k, x_j)_H,
\end{equation}
where $\mathbf E_\ld$ is a spectral function of the operator $A$.
Moreover, the correspondence between all spectral functions of $A$ and all solutions
of the moment problem is one-to-one.
\end{thm}
{\bf Proof. } Only the last statement of the theorem was not proved yet. Its proof repeats the
corresponding proof of Theorem~\ref{t2_2}.
$\Box$

We need some known facts from~\cite{cit_1000_S}.
Let $B$ be a closed symmetric operator in a Hilbert space $H$ with the domain $D(B)$,
which is not necessarily dense in $H$.  Set $\Delta_B(\ld) = (B-\ld E_H) D(B)$,
and $N_\ld = N_\ld (B) = H\ominus \Delta_B(\ld)$, $\ld\in\CC\backslash\R$.

Define an operator $X_i$: $N_i\rightarrow N_{-i}$ in the following way:
\begin{equation}
\label{f3_23}
\varphi = X_i \psi,
\end{equation}
if $\psi\in N_i$, $\varphi\in N_{-i}$ and $\varphi - \psi\in D(B)$.

The operator $X_i$ can be defined also in the following way:
\begin{equation}
\label{f3_24}
D(X_i) = P^H_{N_i} (H\ominus\overline{ D(B) }),
\end{equation}
\begin{equation}
\label{f3_25}
X_i P^H_{N_i} h = P^H_{N_{-i}} h,\qquad h\in H\ominus\overline{ D(B) }.
\end{equation}
The operator $X_i$ is called {\it forbidden with respect to the operator $B$}.

An operator $V$: $N_i \rightarrow N_{-i}$, is called
{\it admissible with respect to the operator $B$}, if
inclusion $V\psi - \psi\in D(B)$ is possible only if $\psi = 0$.
It is equivalent to the condition that the relation $V\psi = X_z \psi$ is possible only if $\psi = 0$.

The formulas
\begin{equation}
\label{f3_26}
D(G) = D(B) \dotplus (V-E_H) D(V),
\end{equation}
\begin{equation}
\label{f3_27}
G( f+V\psi - \psi) = Bf + iV\psi + i\psi,\qquad f\in D(B),\ \psi\in D(V),
\end{equation}
establish a one-to-one correspondence between a set of all admissible with respect to $B$
isometric operators $V$, $D(V)\subseteq N_i$, $R(V)\subseteq N_{-i}$,
and a set of all symmetric  extensions $G$ of the operator $B$.
The operator $G$ is self-adjoint if and only if $D(G)=N_i$, $R(G)=N_{ -i }$.

Denote by $K (B;\CC_+;N_i,N_{-i})$ a class of analytic operator-valued functions $F(\ld)$ in $\CC_+$,
whose values are  contractions which map $N_i$ into $N_{-i}$, $\| F(\ld) \| \leq 1$.

Set $\CC_+^\ve := \{ z\in\CC_+:\  \ve \leq \arg z\leq \pi -\ve \}$, $0\leq\ve\leq \frac{\pi}{2}$.
A function $F\in K (B;\CC_+;N_i,N_{-i})$ is called {\it admissible with respect to the operator $B$}
if relations
\begin{equation}
\label{f3_28}
\lim_{\ld\in \CC_+^\ve,\ \ld\rightarrow\infty} F(\ld)\psi = X_i \psi,
\end{equation}
\begin{equation}
\label{f3_29}
\underline\lim_{\ld\in \CC_+^\ve,\ \ld\rightarrow\infty} \left( |\ld| ( \|\psi\| - \| F(\ld)\psi \| ) \right)
< +\infty,
\end{equation}
imply $\psi = 0$.

The class of all functions from $K (B;\CC_+;N_i,N_{-i})$ which are admissible with respect to
$B$ we denote by $K_a (B;\CC_+;N_i,N_{-i})$.
Notice that in the case $\overline{D(B)}=H$ we have $K_a (B;\CC_+;N_i,N_{-i}) = K (B;\CC_+;N_i,N_{-i})$.

Let $F(\ld)\in K_a (B;\CC_+;N_i,N_{-i})$. In this case the operator $F(\ld)$
is admissible with respect to $B$ \cite{cit_1000_S}. By $B_{F(\ld)}$ we mean an operator $G$ defined as
in~(\ref{f3_26}) with $V=F(\ld)$.

The following theorem holds true, see~\cite[Theorem 12]{cit_1000_S}.
\begin{thm}
\label{t3_2}
Let $B$ be a closed symmetric operator in a Hilbert space $H$ with the domain $D(B)\subseteq H$.
The formula
\begin{equation}
\label{f3_30}
\mathbf R_\ld = \left\{ \begin{array}{cc} (B_{F(\ld)} - \ld E_H)^{-1}, & \img\ld > 0\\
(B_{ F^*(\overline{\ld}) } - \ld E_H)^{-1}, & \img\ld < 0 \end{array}\right.,
\end{equation}
establishes a one-to-one correspondence between the set of all generalized resolvents of $B$ and
the class $K_a (B;\CC_+;N_i,N_{-i})$.
\end{thm}
Using Theorems~\ref{t3_1} and~\ref{t3_2} we get a description of all solutions of the
truncated matrix Hamburger moment problem.
\begin{thm}
\label{t3_3}
Let a truncated matrix Hamburger moment problem~(\ref{f1_5}) with $d\in\N$ be given and
conditions~(\ref{f1_7}) are true. Let the operator $A$ be constructed for the
moment problem as in~(\ref{f3_4}).
All solutions of the truncated moment problem have the following form
\begin{equation}
\label{f3_31}
M(x) = (m_{k,j} (x))_{k,j=0}^{N-1},
\end{equation}
where $m_{k,j}$ satisfy the following relation
\begin{equation}
\label{f3_32}
\int_R \frac{1}{x-\ld} d m_{k,j} (x) = ( (A_{F(\ld)} - \ld E_H)^{-1} x_k, x_j)_H,\qquad \ld\in\CC_+,
\end{equation}
where $F(\ld)\in K_a (A;\CC_+;N_i,N_{-i})$.

On the other hand, to any operator function $F(\ld)\in K_a (A;\CC_+;N_i,N_{-i})$
it corresponds by relation~(\ref{f3_32}) a solution of the truncated matrix Hamburger moment problem.
Moreover, the correspondence between $K_a (A;\CC_+;N_i,N_{-i})$ and
all solutions of the truncated moment problem, established by relation~(\ref{f3_32}), is one-to-one.
\end{thm}
{\bf Proof. } To check the last statement of the theorem it is enough to repeat the arguments
from the proof of Theorem~\ref{t2_4}.
$\Box$

\section{Solvability of the scalar truncated moment problem with even number of given moments.}
Let a moment problem~(\ref{f1_8}) be given. Set
\begin{equation}
\label{f4_1}
\Gamma_n = \left(
\begin{array}{cccc} s_0 & s_1 & \ldots & s_n\\
s_1 & s_2 & \ldots & s_{n+1}\\
\vdots & \vdots & \ddots & \vdots\\
s_n & s_{n+1} & \ldots & s_{2n}\end{array}
\right),\qquad n=0,1,...,d.
\end{equation}
Let $\sigma(x)$ be a solution of the moment problem.
If we choose an arbitrary polynomial $p(x)\in\PP_d$, and calculate
$\int_\R |p(x)|^2 d\sigma(x)\geq 0$,
we can easily see that
\begin{equation}
\label{f4_2}
\Gamma_d \geq 0,
\end{equation}
and therefore all matrices $\Gamma_n$, $0\leq n\leq d$, are real positive semi-definite.
Thus, condition~(\ref{f4_2}) is necessary for the solvability of the moment problem.

Suppose now that a moment problem~(\ref{f1_8}) is given and condition~(\ref{f4_2}) is true.
If $\Gamma_0 = s_0 = 0$, then there exists a unique solution $\sigma(x)=0$, if $s_1=s_2=...=s_{2d+1}=0$,
or there are no solutions in the opposite case.

Assume that $\Gamma_0 = s_0 > 0$. Set
$$ r := \max\{ n:\ 0\leq n\leq d,\ \det\Gamma_n>0 \},\quad 1\leq r\leq d. $$

a) {\it Case $r=d$. } In this case $\Gamma_d > 0$. We can define a real number $s_{2d+2}$ such that
\begin{equation}
\label{f4_3}
\det\Gamma_{d+1} > 0,\quad \Gamma_{d+1} := \left(
\begin{array}{cccc} s_0 & s_1 & \ldots & s_{d+1}\\
s_1 & s_2 & \ldots & s_{d+2}\\
\vdots & \vdots & \ddots & \vdots\\
s_{d+1} & s_{d+2} & \ldots & s_{2d+2}\end{array}
\right).
\end{equation}
To show that, expand the latter determinant by the elements of the last row and choose $s_{2d+2}$ sufficiently
large.
Thus, in this case, by results of V.G.~Ershov and H.~Dym (see the Introduction)
and also by results in~\cite{cit_600_Akh} on the truncated Hamburger moment problem,
it follows that the moment problem~(\ref{f1_8}) has a solution.

b) {\it Case $r<d$. } In this case we have
\begin{equation}
\label{f4_31}
\Gamma_r > 0,\quad \det\Gamma_{r+1} = 0.
\end{equation}
Let $\vec c = (c_0,c_1,\ldots,c_{r+1})$ be a non-zero real vector such that
\begin{equation}
\label{f4_32}
\Gamma_{r+1}\vec c^* = 0,\quad c_{r+1}=1.
\end{equation}
Consider a non-zero real polynomial $p(x) = \sum_{k=0}^{r+1} c_k x^k$,
of degree exactly $r+1$.

If there exists a solution $\sigma(x)$, then
\begin{equation}
\label{f4_5}
\int_\R p^2(x) d\sigma(x) = \sum_{k,n=0}^{r+1} c_k c_n s_{k+n} = 0.
\end{equation}
This implies that $\sigma(x)$ has points of increase only in zeros of $p(x)$,
which we shall denote by $x_0,x_1,\ldots, x_{r}$.
Roots of the polynomial $p(x)$ in this case are real and distinct (or we could replace
$p(x)$ by a polynomial of a less degree such that~(\ref{f4_5}) held, this contradicts~(\ref{f4_31})).
Thus, $\sigma(x)$ is a piecewise constant function, $\sigma(-\infty)=0$,
with jumps in a real distinct points $\{ x_k \}_{k=0}^{r}$.
Denote the jump of $\sigma$ at $x_k$ by $\mu_k$, $0\leq k\leq r$.
The moment equalities~(\ref{f1_8}) are equivalent to
\begin{equation}
\label{f4_6}
\sum_{k=0}^{r} x_k^n \mu_k= s_n,\qquad n=0,1,..., r;
\end{equation}
\begin{equation}
\label{f4_7}
\sum_{k=0}^{r} x_k^n \mu_k= s_n,\qquad n=r+1,r+2,..., 2r+1.
\end{equation}
The linear  system of equations~(\ref{f4_6}) has a non-zero Vandermonde's determinant, and has a unique solution.
This solution should satisfy relations~(\ref{f4_7}).
%From the above considerations we get the following theorem.
\begin{thm}
\label{t4_1}
Let a truncated Hamburger moment problem~(\ref{f1_8}) be given.
It has a solution if and only if

\noindent
a) $s_k = 0$, $k=0,1,...,2d+1$;

\noindent
or

\noindent
b) $\Gamma_d > 0$;

\noindent
or

\noindent
c) $1\leq r<d$, where $r := \max\{ n:\ 0\leq n\leq d,\ \det\Gamma_n>0 \}$;
the polynomial $p(x) = \sum_{k=0}^{r+1} c_k x^k$, where $c_k$ are complex numbers satisfying~(\ref{f4_32}),
has real distinct zeros $\{ x_k \}_{k=0}^r$; and the unique solution of linear system~(\ref{f4_6})
consists of non-negative numbers $\mu_k\geq 0$, $k=0,1,...,r$, which satisfy relations~(\ref{f4_7}).

In cases a) and c) the solution is unique.
\end{thm}
{\bf Proof. } The necessity follows from the above considerations. The sufficiency of condition~b) was shown.
The sufficiency of conditions~a) and c) is obvious.
$\Box$

\begin{center}
\large\bf
A description of all solutions of the matrix Hamburger moment problem in a general case.
\end{center}
\begin{center}
\bf S.M. Zagorodnyuk
\end{center}

We describe all solutions of the matrix Hamburger moment problem in a general case (no conditions
besides solvability are assumed). We use the fundamental results of A.V.~Shtraus on the generalized
resolvents of symmetric operators. All solutions of the truncated matrix Hamburger moment problem with an odd number of
given moments are described in an "almost nondegenerate" case. Some conditions of solvability
for the scalar truncated  Hamburger moment problem with an even number of given moments are given.

MSC: 44A60; Secondary 30E05
Key words:  moment problem, positive definite kernel, spectral function.

\end{document}